# Truth-relevant Logic – Propositional Calculus

*X.Y. Newberry*

## Abstract

The thesis of this paper is that truth-relevant logic is a better foundation for mathematics than classical logic. It is a system  proposed by Richard Diaz in 1981. In a certain sense t-relevant logic is based on Kleene's strong tables. These define a system with three values: true, false, unknown. It turns out that there exist tautologies with the following property: there exists a proper subset of propositional variables (t-relevant variables) such that for all combinations true/false the tautology will be true, that is, the rest of the variables (redundant variables) occurring in the tautology can be unknown. We consider such compound sentences as neither true nor false. Philosophical justification is provided. Proof system based on tableaux is proposed. The following theorem is proved: tautology $L = (R_1 \lor \sim R_1) \ \& \ (R_2 \lor \sim R_2) \ \& \ ... \ \& \ (R_m \lor \sim R_m)$, where $\{R_1, R_2, ... R_m\}$ is a subset of t-relevant propositional variables, and '=' is equivalence according to Kleene's strong tables.

# 1. Introduction

Consider the formula

$$\sim(\exists x)(\exists z)(Prf(x,\ z)\ \&\ Diag(y,\ z)) \tag{U}$$

with one free variable *y*,, where *Diag(y, z)* means that *y* is the Gödel number of a formula with one free variable iff *z* is the Gödel number number of the formula obtained from *y* by substituting (the numeral of) the Gödel number of *y* for the free variable in *y*. And *Prf(x, z)* means that *x* is the Gödel number of a sequence that is a proof of the sentence with Gödel number *z*.

Let the constant *k* be the Gödel number of (U). We substitute *k* for the free variable *y* in (U), and obtain

$$\sim(\exists x)(\exists z)(Prf(x,z)\ \&\ Diag(k,z)) \tag{G}$$

This is Gödel's self-referential sentence. Let $\ulcorner G \urcorner$ be the Gödel number of (G). In classical logic

$$\sim(\exists x)Prf(x,\ \ulcorner G \urcorner) \tag{H'}$$

is equivalent to (G). But I argued in Newberry (2015) that in Strawson's logic of presuppositions

$$(\exists x)Prf(x,\ \ulcorner G \urcorner) \tag{H}$$

is a *presupposition* of (G). Hence when (H) is false then (G) is *neither true nor false*.

The question is how to formalize a logic where we could derive (H') but not (G). A key observation is that Gödel's sentence (G) is *vacuous*. We therefore need a logic where vacuous sentences are not derivable. In such a logic Gödel's self-referential sentence will not be derivable, but *∼(∃x)Prf(x,$\ulcorner G \urcorner$)* might be.

We need to start at the propositional level. And it turns out that a logic suitable for our purposes already exists. In 1981 M. Richard Diaz published a monograph titled *Topics in The Logic of Relevance*. He proposed a new type of relevance logic called *truth-relevant logic*. In this logic the paradox of material implication *(P & ∼P) → Q* is *not* truth-relevant, but neither are its classical equivalents such as *(P ∨ ~P) ∨ Q*. By generalization of the former we obtain *(∀x)(Px & ~Px) → Qx)*, i.e. a *vacuous* sentence which ought not to be derivable. Newberry (2019b) defines the semantics of such predicate logic while (2019c) outlines a derivation system

Diaz conceived his system as a type of a relevance logic. We will look at it from a different angle. In section 3 I attempt a philosophical justification of t-relevant logic. Diaz (pp. 74-75.) described a derivation system, called *SR*. I do not find it particularly elegant or intuitive. Nor is it immediately obvious how it could be extended to the predicate calculus. In section 4 I propose a more streamlined system based on truth trees.

Diaz did not have much faith in this system; his investigations have a hypothetical character, perhaps something to ponder for further inspiration. I do not share this skepticism. I believe that truth-relevant logic can become the new foundation of logic and mathematics. One objection he raised is addressed in Appendix A.

This paper is *not* about mathematics but about the *foundations of* mathematics. The philosophy is discussed in Newberry (2014). It probably should have been made more explicit that it is based on l*ogical positivism*, perhaps with some modifications and qualifications, in particular on (unqualified) l*ogicism* and profound skepticism regarding the existence of any mathematical Platonic entities.

## 2. Preliminaries

Here is a brief recapitulation of the basic *definitions* according to Diaz:

**Definition 1:** A set of propositional variables is *truth determining* for a proposition G iff the truth value of G may be determined as true or false on all assignments of T, F to the set.

**Definition 2:** Let $p_1, p_2, \ldots p_n$ be all the variables occurring in G, $p_i$ is *truth-redundant* in G iff there is a truth determining set for G that does not contain $p_i$. G is *truth-relevant* iff G contains no truth-redundant variables.

Of special interest is the set of tautologies that are also t-relevant.

(Diaz, 1981, pp. 66, 67)

The key idea is this. Given $A \lor B$, if $A$ is true we know that the formula is true *regardless* of the value of $B$. Given $A \& B$, if $A$ is false we know that that the formula is false *regardless* of the value of $B$.

Consider $P \rightarrow (Q \rightarrow P)$. If $T(P)$ then $Q \rightarrow P$ is true regardless of the value of $Q$. But then the entire expression is true. If $F(P)$ then the expression is true regardless of the value of $(Q \rightarrow P)$. If $T(Q)$ then we cannot determine the truth value of the formula without knowing the value of $P$. Thus $\{P\}$ is truth determining, $Q$ is truth-redundant. The formula is *not* t-relevant. (Diaz, p.65.)

In truth relevant logic the truth tables are constructed from the following elementary truth tables (Diaz, p.66.)

```
V | T x F      & | T x F       → | T x F      ~ |
--+-------     --+-------      --+-------     ---+---
T | T T T      T | T x F       T | T x F      T | F
x | T x x      x | x x F       x | T x x      x | x
F | T x F      F | F F F       F | T T T      F | T
```

                    Table 2.1

These tables happen to be identical with Kleene's "strong" tables. But Diaz's objective is to identify irrelevant variables (sentence letters, atomic formulas), i.e. variables that do not make any contribution to the truth value of the compound sentence.

Here is the truth table for $((A \rightarrow B) \rightarrow A) \rightarrow A$:

```
   A | B | ((A -> B) -> A) -> A
  ---+---+----------------------
1) T | T |      T        T      T
2) T | F |      F        T      T
3) F | T |      T        F      T
4) F | F |      T        F      T
```

                    Table 2.2

Below is the truth table when the value of $B$ is unknown. It was obtained in accordance with the elementary tables in Table 2.1.

```
     A | B | ((A -> B) -> A) -> A
   --+---+----------------------
1) T | x |       x       T    T
3) F | x |       T       F    T
```

Table 2.3

It is apparent that the value of *B* is not required to determine that the formula is a tautology. The same cannot be said about *A* as the next table shows.

```
     A | B | ((A -> B) -> A) -> A
   --+---+----------------------
1) x | T |       T       x    x
2) x | F |       x       x    x
```

Table 2.4

Thus {*A*} is truth determining, *B* is truth redundant, and *((A → B) → A) → A* is not t-relevant.

It is illustrative to picture the formula as a switching circuit.

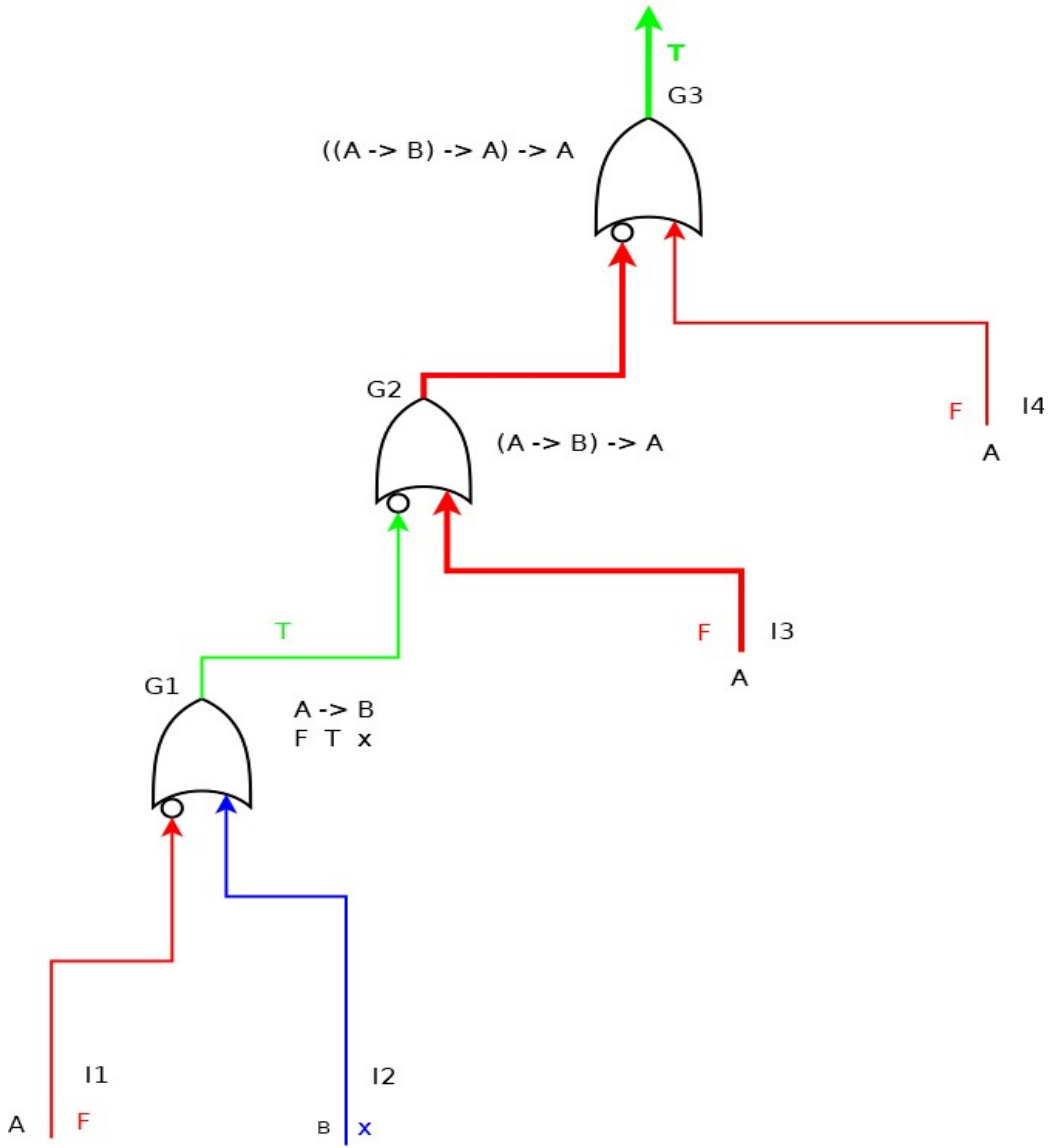

Figure 2.1

Note that if *F*(*A*) then the *B* signal cannot propagate past the first gate. Figure 2.1. A definite truth value appears under the first implication in Table 2.3.

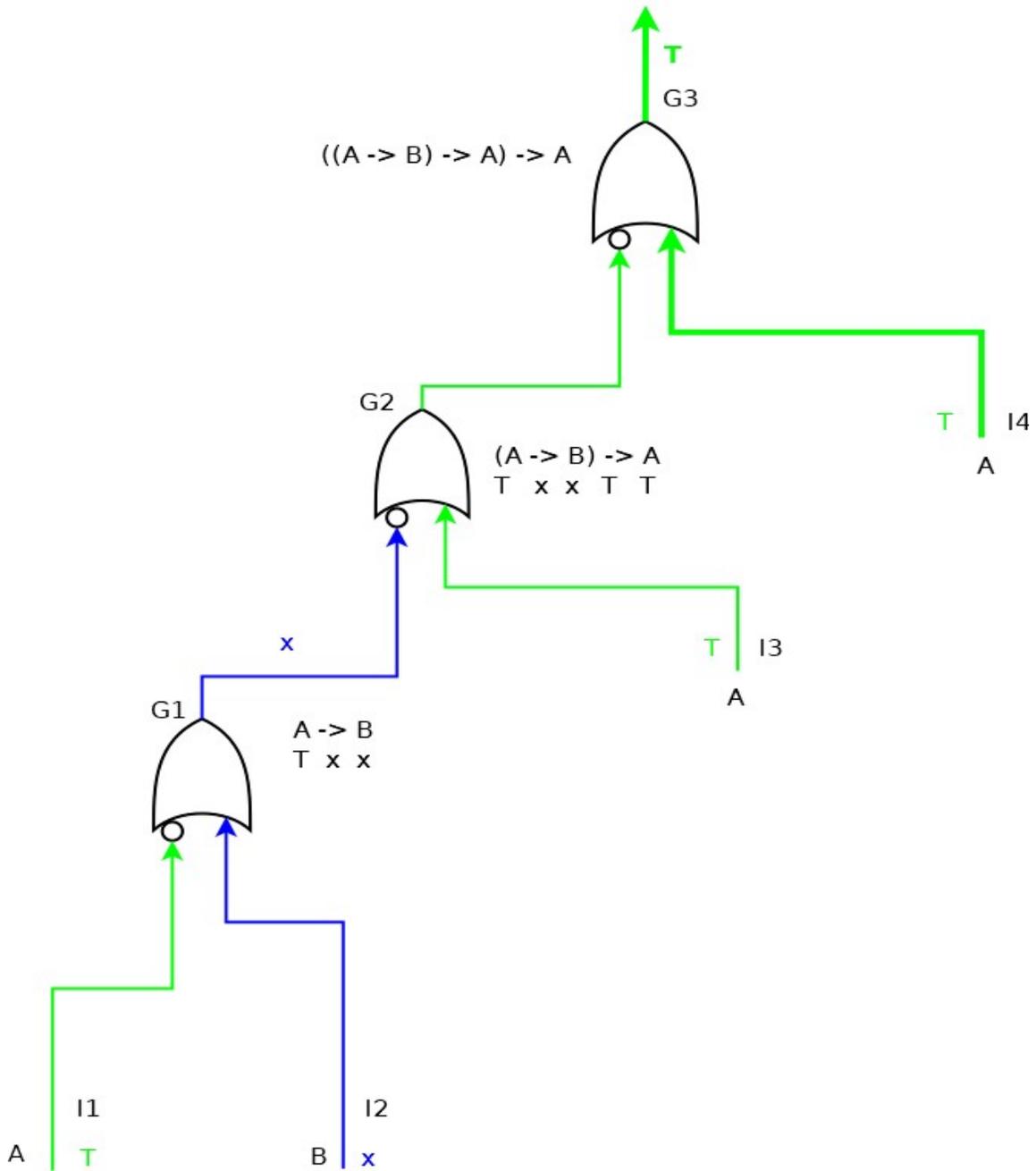

Figure 2.2

If $T(A)$ then the $B$ signal cannot propagate past the second gate. Figure 2.2. A definite truth value appears under the second implication in Table 2.3. That is as far as the $B$ signal (represented by $x$) gets. The path from the input $B$ cannot be sensitized. Under no conditions will the $B$ signal influence the output.

# 3. Justification/Interpretation

Classical logic views '∨' and '&' as purely Boolean operators operating on the truth values of propositions. This actually may not be the best approximation of 'or' and 'and' as used in the natural language. I would like to suggest a somewhat modified interpretation. The expression $P \lor Q$ may be better understood as indicating the speaker's belief that it is *possible* that *F(P)* and that it is *possible* that *F(Q)*, but at least one of them is true. Indeed, if the speaker were sure that *T(P)* why would he append '∨ *Q*'? It is is not necessary to construct any formal or rigorous theory of modality, i.e. a theory of what *possible* precisely means. It simply indicates the speaker's belief just as uttering *P* indicates his belief that *P* is true without a necessity to elaborate how *P* may be verified.

Analogically I would like to suggest that *~(P & Q)* ought to be understood as the speaker's belief that it is *possible* that *T(P)* and that it is *possible* that *T(Q),* but at least one of them is false. Similar considerations apply to implication $P \to Q$.

Nobody in his right mind would utter the sentence "*All bachelors are unmarried or it is raining outside*." Only a logician is capable of something like that. But there is something much more fundamental going than merely proximity to the natural language. Consider

$$(P \& \sim P) \to Q \tag{3.1}$$

This is the *paradox of material implication* – a counterintuitive result although a tautology of classical logic. The sentence is attempting to say what would happen if the impossible happened. Such a situation cannot be *pictured*, and hence the above expression ought to be considered *meaningless*. We simply cannot form a *mental image* of what '⊥ → *Q*' is trying to convey..

One could object that *P & ~P* cannot be pictured either. But that is precisely how *~(P & ~P)* should be read, i.e. as an asserting  "*(P & ~P) cannot be pictured*." Then (3.1) would be attempting to say something to the effect "*When the circumstances that cannot be pictured occur then P*." It can be said that something is impossible, but it cannot not be said what would happen if the impossible happened. It cannot be said because we cannot form a *mental image* of it.

But (3.1) is equivalent to

$$(P \lor \sim P) \lor Q \qquad\qquad\qquad (3.2)$$

It says something to the effect: "*The situation we cannot picture <u>otherwise</u> is the case or it is raining outside.*" It does not strike us as absurd as (3.1) although it is still difficult to form a *mental image* of what exactly we are trying to say here. Nevertheless the flip side is still (3.1), i.e. "*If the situation we cannot picture is the case then it is raining outside.*" Therefore we ought to consider (3.2) meaningless as well.

<p style="text-align:center">* * * * * * *</p>

Suppose you read in a scientific magazine "*There are amino acids on Proxima Centauri b.*" You would be surprised. (The statement carries a lot of information.) Now suppose you read in a scientific magazine "*Square is a rectangle or there are amino acids on Proxima Centauri b.*" You would be even more surprised! The author apparently is trying to convey that the compound sentence is true. But we know the first simple sentence is true, and therefore we know that the compound sentence is true regardless of the truth value of the second component. The compound sentence *conveys no information* about amino acids on Proxima Centauri b. The second component is meaningless in this sense. That is, $G$ carries a lot of information on its own but zero information when a part of $\top \lor G$. So in fact we do not have $T \lor T$ (or $T \lor F$), we have $T \lor M$, where $M$ stands for meaningless. It is reasonable to put $V * M = M$, where $V$ is any truth value and '*' is any connective. So we have the result $M(\top \lor G)$; $\top$ radiates dangerous X-rays, which bleach everything in their path.

<p style="text-align:center">* * * * * * *</p>

I am proposing the following interpretation.

$$|P \text{ v } Q| = T \ \text{ iff } \ (P + Q) \bullet \Diamond\sim P \bullet \Diamond\sim Q \qquad\qquad (3.3)$$

$$|P \text{ v } Q| = F \ \text{ iff } \ \sim(P + Q) \bullet \Diamond\sim P \bullet \Diamond\sim Q \qquad\qquad (3.4)$$

where '+', '•' are a boolean operators, '$\Diamond$' means *possible*. In case of contingent sentences possibility is just the belief of the speaker. |$A$| indicates the truth value of *A*. When $\Box P + \Box Q$ then there is a truth value gap. (The symbol '$\Box$' means *necessary*.)

We will call "◇~P • ◇~Q" the *conditions of relevance*. They play a similar role as presuppositions in Strawson's logic.

But *~(P + Q) = ~P • ~Q*. So in case of (3.2) the conditions of relevance are satisfied automatically; *P v Q* cannot be false unless they are satisfied. The clause ◇~P • ◇~ Q is redundant and may be omitted.

Furthermore we will define"

|P & Q| = T  iff  P • Q • ◇P • ◇Q                                    (3.5)

|P & Q| = F  iff  ~(P • Q) • ◇P • ◇Q                                 (3.6)

But again, in the first case the conditions of relevance are satisfied automatically; *P & Q* cannot be true unless they are satisfied.  Furthermore we will define for any formula S:

|S| = T  iff |~S| = F                                                (3.7)

|~S| = T  iff |S| = F                                                (3.8)

For example:

|~(P & Q)| = T  iff  |P & Q| = F, i.e.                               (3.9)

|~(P & Q)| = T  iff  ◇P • ◇Q • ~(P • Q)                              (3.10)

We will summarize our definitions in a table.

| |S| = T | Boolean function | Relevance conditions | |S| = F | Boolean function | Relevance conditions |
|---|---|---|---|---|---|---|
| **|P v Q| = T** | P + Q | ◇~P • ◇~Q | **|P v Q| = F** | ~(P + Q) | Automatic |
| **|P & Q| = T** | P • Q | Automatic | **|P & Q| = F** | ~(P • Q) | ◇P • ◇Q |
| | | | | | |
| |~(P v Q)| = T | ~(P + Q) | Automatic | |~(P v Q)| = F | P + Q | ◇~P • ◇~Q |
| |~(P & Q)| = T | ~(P • Q) | ◇P • ◇Q | |~(P & Q)| = F | P • Q | Automatic |

Table 3.1

The somewhat vague concept of what is *possible* acquires rather a crisp meaning when we talk about Boolean functions. Let us suppose that any input can be either 0 or 1. But given $C = A + {\sim}A$ and $D = C + B$ then *C = 0* is *impossible*. The same applies to compound formulas such as

$$(P \bigvee {\sim}P) \bigvee Q \tag{3.11}$$

Clearly $\Box(P \vee {\sim}P)$. Therefore the condition $\Diamond{\sim}(P \vee {\sim}P)$ is not satisfied. Hence (3.11) is not true. When building formal theories we will assume that all contingent sentences are possible. We will define

$$\Diamond G \ \ \text{=def=} \ \ \nvDash{\sim}G, \ \ \Box G \ \ \text{=def=} \ \ \vDash G \tag{3.12}$$

* * * * *

We have observed that t-relevant logic gets rid of the prototypical version of the *paradox of material implication* (3.1). Let us now look at a more complicated paradox from Wikipedia:

> If Nadia is in Barcelona then she is in Madrid or if Nadia is in Madrid then she is in Barcelona.

Needless to say, t-relevant logic makes a short work of it. It can be formalized as

$$({\sim}B \bigvee M) \bigvee ({\sim}M \bigvee B)$$

It *is* a tautology, but *not* a t-relevant tautology as apparent from the table below.

```
B M | ( ~ B ∨ M ) ∨ ( ~ M ∨ B )
----+------------------------------
T x |   F T x x   T   x x T T
F x |   T F T x   T   x x x F
```
<div align="center">Table 3.1</div>

So much for that. But what else can we say about the offending sentence? Why does the argument seem absurd to us even though it is a tautology according to classical logic? Because Madrid and Barcelona do not overlap! According to my "modal" interpretation:

$${\sim}B \bigvee M \ = \ \Diamond \boldsymbol{B} \bullet \Diamond{\sim}M \bullet ({\sim}B + \boldsymbol{M})$$

where '+', '•' are purely Boolean operators operating on truth values. That is, the sentence asserts (among other things) that Nadia is in Madrid, and that it is possible that at the same time she is in Barcelona. We know that it is not. We perceive the original sentence as paradoxical because it evokes the image of Nadia in Barcelona. Whether Nadia is actually in Barcelona or not we interpret the sentence hypothetically: If Nadia were in Barcelona then she could not be in Madrid. We can picture her in Barcelona, but *not* simultaneously in Madrid.

* * * * * * *

One potential objection to t-relevant logic, suggested by its own originator, is that given *A & B* we cannot extract either *A* or *B*, as *A & B → A* is not a t-relevant tautology. That is, the expression *A & B* cannot be "dissolved". But if *T(A & B)* then certainly *T(A)* and *T(B)*. The solution of this dilemma probably is that

$$A \& B \vdash A, B$$

is not the same as

$$A \& B \to A, \; A \& B \to B$$

This implies that unrestricted conditional proof cannot be used in t-relevant logic. [1])

---

# 4. A Proof System – Propositional Calculus

Let us now turn our attention to the problem of constructing a proof system for t-relevant logic. One such system, called *SR*, can be found in Diaz (pp. 74-75.) As mentioned above, this system has some drawbacks. I would like to propose a more streamlined method based on truth trees. The tree method is described in general e.g. in Bostock (1997). For the time being we will focus on proving tautologies rather than sequents generally. We will use the term *tableau* and *truth tree* interchangeably. The key observation is that the tableau method, as we know it, already distinguishes between t-relevant and non t-relevant tautologies. It shows that the negation of the compound formula *L* to be proven implies *($R_1$ & ~$R_1$) ∨ ($R_2$ & ~$R_2$) ∨ ... ∨ ($R_m$ & ~$R_m$)*, where the $R_i$ are a subset *S* of the variables (sentence letters) occurring in *L. If S is a **proper subset** then then the formula L is **not** t-relevant.* We will further elaborate on this observation.

**Theorem 4.1** Let $L(P_1, P_2, ... P_n)$ be any tautology such that $P_1, P_2, ... P_n$ are all the variables occurring in *L*, and let our connectives be *interpreted according to the tables 2.1*, i.e. *x* is a third value. Then *L* is equivalent to

$$(R_1 ∨ \sim R_1) \& (R_2 ∨ \sim R_2) \& ... \& (R_m ∨ \sim R_m) \tag{4.1}$$

if and only if $\{R_i\}$ is a set of t-relevant variables occurring in *L*. Naturally $\{R_i\} \subseteq \{P_i\}$.

**Proof:** First we prove that $\{R_i\}$ contains all the t-relevant variables. So suppose that there exists a t-relevant variable *G* occurring in *L* but not in $\{R_i\}$. Suppose further that for all $R_i$, $|R_i| = T$ or $|R_i| = F$, and $|G| = x$. Then the truth value of (4.1) is *T* but the truth value of *L* is *x*. Now suppose that $\{R_i\}$ contains a redundant variable *H*, and $|H| = x$, all other variables in $\{R_i\}$ are *T* or *F*. Then the truth value of (4.1) is *x*, but the truth value of *L* is *T*.

```
   Assignment      | 4.1 | L
   -----------------------------
1) G ∉ {Ri}, G = x |  T  | x       G is relevant
2) H ∈ {Ri}, H = x |  x  | T       H is redundant
```

Table 4.1

When $\{R_i\}$ are precisely the variables of a truth determining set then *L* and 4.1 agree as follows.

```
   Assignment | 4.1 | L
   ----------------------------
1)   T v F    |  T  | T
2)     x      |  x  | x
```

<div align="center">Table 4.2</div>

Line 1 means that definite values have been assigned to all the variables in $\{R_i\}$, line 2 means that $x$ has been assigned to at least one $R_i$. QED.

So we have

$$L(P_1, P_2, \ldots P_n) \;=\; (R_1 \bigvee \sim R_1) \;\&\; (R_2 \bigvee \sim R_2) \;\&\; \ldots \;\&\; (R_m \bigvee \sim R_m) \tag{4.2}$$

as well as

$$\sim L(P_1, P_2, \ldots P_n) \;=\; (R_1 \;\&\; \sim R_1) \bigvee (R_2 \;\&\; \sim R_2) \bigvee \ldots \bigvee (R_m \;\&\; \sim R_m) \tag{4.3}$$

The tableau method works as follows.

$$\sim L(P_1, P_2, \ldots P_n) \vdash (Q_1 \;\&\; \sim Q_1) \bigvee (Q_2 \;\&\; \sim Q_2) \bigvee \ldots \bigvee (Q_m \;\&\; \sim Q_m) \tag{4.4}$$

We have to prove that the set $\{Q_1, Q_2, \ldots Q_m\}$ produced by the tableau is a set of t-relevant variables $\{R_1, R_2, \ldots R_n\}$ of $L$. Note that the tree preserves equivalence according to the tables 2.1. For example here are the rules for v and &.

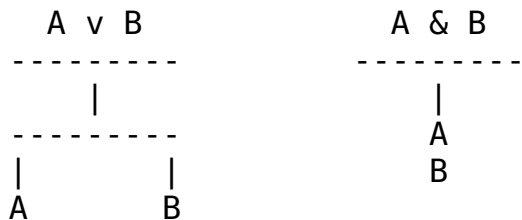

Here $A$ v $B$ is true if $|A| = T$ regardless if $|B| = T, F$ or $x$, $A$ & $B$ is false if $A$ is false regardless if $|B| = T, F$ or $x$ etc. In fact the tree is but the original formula in a graphically different form, which makes it easier to spot self-contradictions such as $Q_i \;\&\; \sim Q_i$. Only instead of v there is branching and instead of & there is vertical juxtaposition. When the tree is closed the result is

$$(Q_1 \;\&\; \sim Q_1 \;\&\; E_1) \bigvee (Q_2 \;\&\; \sim Q_2 \;\&\; E_2) \bigvee \ldots \bigvee (Q_m \;\&\; \sim Q_m \;\&\; E_m) \tag{4.5}$$

where the $E_i$ is a conjunction of the expressions on the last edge of the *i*-th branch.. One disjunct – one branch. The formula above is clearly equivalent to the original formula according to the tables 2.1, and it is further equivalent to

$$(Q_1 \mathbin{\&} \sim Q_1) \bigvee (Q_2 \bigvee \mathbin{\&} \sim Q_2) \bigvee \ldots \bigvee (Q_m \mathbin{\&} \sim Q_m) \qquad (4.6)$$

as $Q_i \mathbin{\&} \sim Q_i$, is false, and it does not matter if $E_i$ is *T, F* or *x*. Since the formula above is equivalent to the original formula then according to Theorem 4.1 $\{Q_i\}$ is a set of t-relevant variables occurring in *L*.

<div align="center">* * * * * * *</div>

Let us look at some examples. Consider $(P \to Q) \to ((Q \to R) \to (P \to R))$ :

```
   P | Q | R |  (P → Q) → ((Q → R) → (P → R)
  ---+---+---+--------------------------------
1) T | T | T |    T    T    T    T    T
2) T | T | F |    T    T    F    T    F
3) T | F | T |    F    T    T    T    T
4) T | F | F |    F    T    T    F    F
5) F | T | T |    T    T    T    T    T
6) F | T | F |    T    T    F    T    T
7) F | F | T |    T    T    T    T    T
8) F | F | F |    T    T    T    T    T
```

<div align="center">Table 4.3</div>

This tautology *is t-relevant* as we can observe from the tables below. When the value of any of the three variables is unknown then the truth value of the formula cannot be determined.

```
   P | Q | R |  (P → Q) → ((Q → R) → (P → R))
  ---+---+---+--------------------------------
1) T | T | x |    T    x    x    x    x
3) T | F | x |    F    T    T    T    x
5) F | T | x |    T    T    x    T    T
7) F | F | x |    T    T    T    T    T
```

<div align="center">Table 4.4</div>

```
    P | Q | R |   (P → Q) → ((Q → R) → (P → R))
 ---+---+---+-------------------------------------
 1) T | x | T |      x     T     T     T     T
 2) T | x | F |      x     x     x     x     F
 5) F | x | T |      T     F     T     F     T
 6) F | x | F |      T     T     x     T     T
```

Table 4.5

```
    P | Q | R |   (P → Q) → ((Q → R) → (P → R))
 ---+---+---+-------------------------------------
 1) x | T | T |      T     T     T     T     T
 2) x | T | F |      T     T     F     T     x
 3) x | F | T |      x     T     T     T     T
 4) x | F | F |      x     x     T     x     x
```

Table 4.6

Note that according to the tables 2.1

$$(P \to Q) \to ((Q \to R) \to (P \to R)) \neq (P \lor {\sim}P) \mathbin{\&} (Q \lor {\sim}Q) \tag{4.7}$$

The assignment $|P| = T$, $|Q| = T$, $|R| = x$ makes the left side $x$ (table 4.4 above) and the right side true. This is a case of $G \notin \{R_i\}$ (table 4.1, line 1.) By negating (4.7) we obtain

$${\sim}((P \to Q) \to ((Q \to R) \to (P \to R))) \neq (P \mathbin{\&} {\sim}P) \lor (Q \mathbin{\&} {\sim}Q) \tag{4.8}$$

And since the tree preserves equivalence according to the tables 2.

$${\sim}((P \to Q) \to ((Q \to R) \to (P \to R))) \nvdash (P \mathbin{\&} {\sim}P) \lor (Q \mathbin{\&} {\sim}Q) \tag{4.9}$$

but rather

$${\sim}((P \to Q) \to ((Q \to R) \to (P \to R))) \vdash (P \mathbin{\&} {\sim}P) \lor (Q \mathbin{\&} {\sim}Q) \lor (R \mathbin{\&} {\sim}R) \tag{4.10}$$

Here is the actual tree:

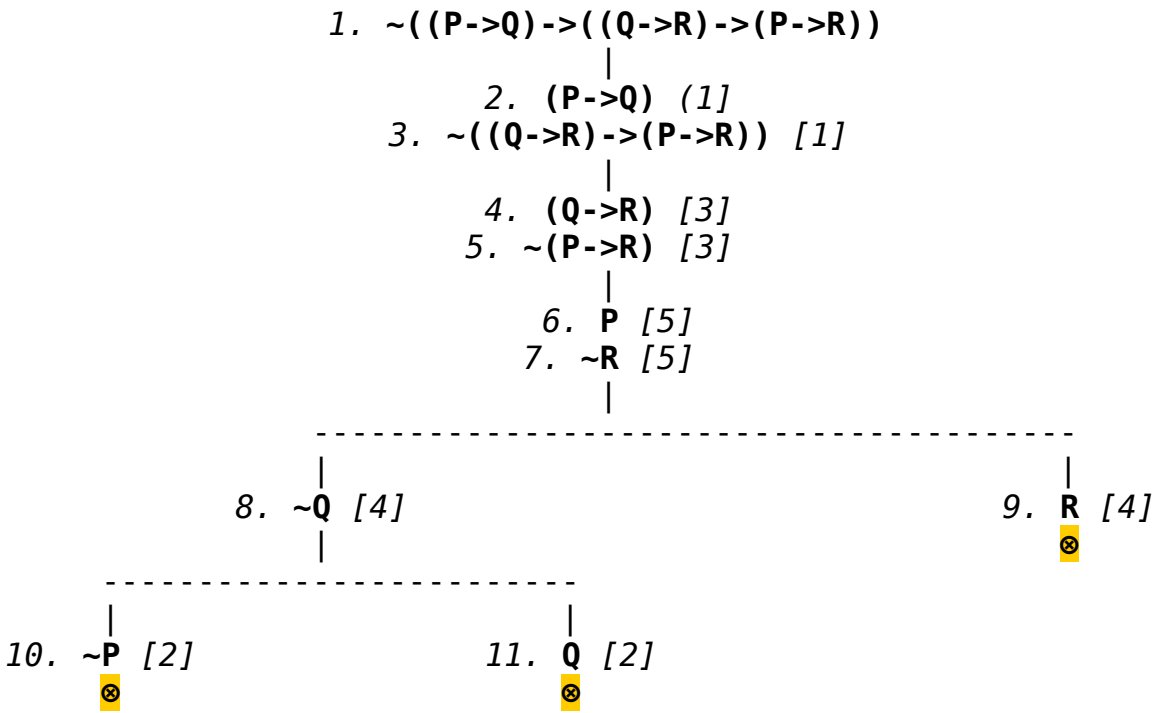

Figure 4.1

\* \* \* \* \* \* \*

Now consider Peirce's law $((A \to B) \to A) \to A$:

```
  A | B | ((A -> B) -> A) -> A
 ---+---+-----------------------
1) T | T |     T      T     T
2) T | F |     F      T     T
3) F | T |     T      F     T
4) F | F |     T      F     T
```

Table 4.7

*This tautology is **not** t-relevant.* (Peirce's law is not a theorem of intuitionistic logic either, there are some similarities.)

```
  A | B | ((A -> B) -> A) -> A
 ---+---+-----------------------
1) T | x |     x      T     T
3) F | x |     T      F     T
```

Table 4.8

```
    A | B | ((A -> B) -> A) -> A
  ---+---+----------------------
1) x | T |     T      x      x
2) x | F |     x      x      x
```

<div align="center">Table 4.9</div>

In particular $B$ is truth-redundant, and $\{A\}$ is the truth determining subset. Note that

<p style="color:red">$((A \to B) \to A) \to A \;\neq\; (A \bigvee \sim A)\ \&\ (B \bigvee \sim B)$</p>

$$(4.11)$$

The assignment $|A| = T$, $|B| = x$ makes the left side true (table 4.8) but the right side $x$
This is a case of a redundant $H \in \{R_i\}$, table 4.1, line 2. The following is correct.

$$((A \to B) \to A) \to A \;=\; A \bigvee \sim A \tag{4.12}$$

And since the tree preserves equivalence according to the tables 2.1

$$((A \to B) \to A) \to A \;\vdash\; A \bigvee \sim A \tag{4.13}$$

Here is the actual tree:

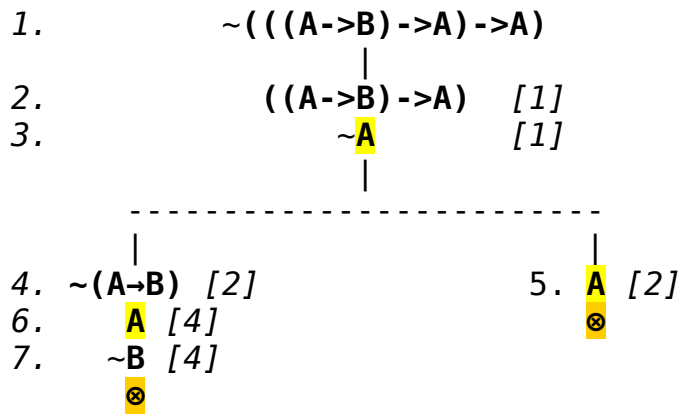

<div align="center">Figure 4.2</div>

We have obtained $A\ \&\ \sim A\ \&\ B$. This is equivalent to $A\ \&\ \sim A$; $B$ can be whatever it wants, even $x$. This is how the tree eliminates irrelevant variables. Here $B$ plays the role of $E_1$ in formula (4.5) mentioned above.

<div align="center">* * * * * * *</div>

A formula can have multiple distinct truth-determining subsets. A case in point is
*(P ∨ ~P) ∨ (Q ∨ ~Q)*. It is *not* t-relevant, and it has *two* truth-determining disjoint proper
subsets according to definition 1, namely *{P}* and *{Q}*. The same applies to *(~B ∨ M) ∨*
*(~M ∨ B)*, which we considered in section 3. Let us look at the truth tree for the former:

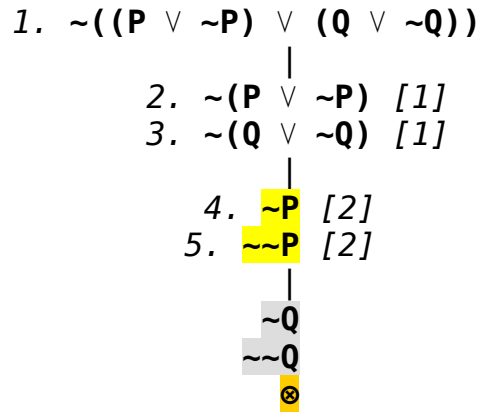

Figure 4.3

We observe that we have obtained *P & ~P & E*. The tree closed before *Q* had a chance to
make any contribution. By convention we close a branch as soon as we encounter the first
contradiction. Here *E = Q & ~Q*, but it is of no further interest to us. We know that
*(P & ~P & E) = (Q & ~Q)* according to the tables 2.1 regardless of what *E* is. Again this is
how the tree eliminates irrelevant variables. So *P* is relevant. The only other variable
occurring in the formula is *Q*; it is redundant.

But we could have as well produced the following tree:

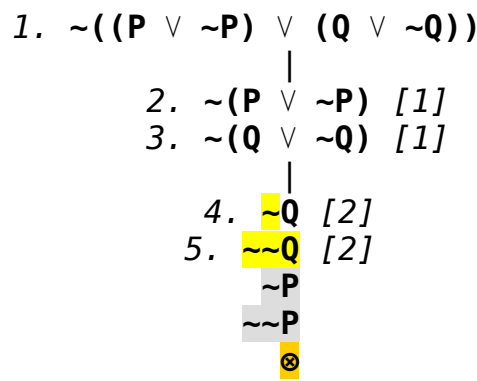

Figure 4.4

We have obtained  *Q & ~Q & E'*. In this case *Q* is relevant and *P* redundant. There are two ways about it but *P* and *Q* are not both relevant at the same time. The only way a variable *V* can be relevant is when *V* ∨ ~*V* appears on a branch with *no other* contradiction. When *V* appears only on branches with existing contradictions it will be a) redundant and b) eliminated. For we cherry-pick contradictions and by convention we close a branch as soon as we find the first one. (Although we can change the order, and thus obtain different sets of t-relevant variables.) This is why (4.6) does not contain any redundant variables. Another way of putting it is that (*P & ~P) &  (Q & ~Q)* does not conform to (4.1) while *P & ~P* and  *Q & ~Q* do.

Truth-relevance has been staring at us from the tableaux all along. When all variables in a formula *~L* are self-contradicted we have proven *T(L)*. If the tree closes but not all variables are contradicted we have merely proven *~F(L)*.

# APPENDIX A

We will address an objection Diaz raised against his own system.

> we could not claim to have proved an entailment until we were sure that each of the premises in the theorem were necessary to reach the conclusion. And for that matter, if someone proved a theorem "A and B entails C," and it was found later that an assumption weaker than B, in conjunction with A. was sufficient to derive C, we would have to claim that it is false [2]) that "A and B entails C." . . . The problems may be illustrated by a basic theorem of Group Theory. Suppose we offered the following theorem and proof:

> **Theorem:** If G is a group then the identity element e of G is unique.

> **Proof:** Assume that we have a group G with identity elements e and f. By the properties of a group e=ef=fe and f=fe=ef. Hence e=f and the identity is unique.

> By the relevance conditions of OR [3]), we cannot claim that the theorem is true until we ascertain that all of the premises are relevant to the derivation. But "G is a group" is really a complex statement; one whose real force can only be seen by replacing "group" with its definition. G Is a non-empty set together with a binary operation on this set, such that for all a, b, c

> 1. $ab \in G$
> 2. $a(bc) = (ab)c$
> 3. There exists an $e \in G$ such that $ae=ea=a$ for all $a \in G$.
> 4. For all $a \in G$ there exists an $a^{-1} \in G$ such $aa^{-1} = a^{-1}a = e$,

> Thus the statement "G is a group" is actually a set of statements, which includes statements about the existence of elements of G and a binary operation on those elements, associativity of the operator, closure under that operator, and the existence of two sided identity and inverse elements. But not all of these properties are required for the proof that the identity element is unique; in particular, the existence of inverses is not required for the proof. Hence, it would seem that, since the theorem implicitly contains irrelevant premises, our theorem is incorrect as it stands.

> Now I believe that this example suffices to show that strict adherence to principles of relevance would lead to disaster in logic and mathematics. The requirements for proof would be such as to actually block the development of interesting and useful results instead of promoting them. (Diaz, 1981, pp. 134-135)

OK, needless to say we will not actually try to prove the uniqueness of *e* from all four axioms. We will use only the three relevant ones. Generally we will develop a theory for each subset of the axioms, and we will call the union of these theories the theory of groups.

---

[2] In our conception it would rather be not true.
[3] 'OR', 'occurrence relevance logic' ia a variant of t-relevant logic. (Diaz, 1981, pp. 116 – 126.)

It is more like saying: "Here is a bunch B of axioms, pick the ones you need, and prove your theorem." Once we prove our theorem from the relevant premises we can loosely speaking say that we have proved it from the bunch B. We have to be precise when it is necessary, but we can afford to be lax when it is not. An extra premise is not a disaster – we will simply drop it.